\input AHTOHFIE.STY
\hfuzz11pt


\UDC{512.543.1+512.543.7+512.54.05}
\MSC{20e06,20f05,20f06,20f67,20f70}

\title{
Relative hyperbolicity and similar properties
\\
of
one-generator one-relator relative presentations
with powered unimodular relator
}
\author{%
Anton A. Klyachko
\quad
Denis E. Lurye
}

\address{
\myAddress
\quad
doomden1990@yahoo.com
}
\grants{\RFBR11-01-00945}

\abstract{%
A group obtained from a nontrivial group
by adding one generator and one relator which is a
proper power of a word in which the exponent-sum of the
additional generator is one
contains the free square of the initial group and
almost always (with one obvious exception) contains
a non-abelian free subgroup.
If the initial group is involution-free or
the relator is at least third power,
then the obtained group is
SQ-universal and
relatively hyperbolic
with respect to the initial group.
}

\s 1.
Introduction

Let $G$ be a torsion-free group and let a group $\^G$ be obtained
from the group $G$ by
adding one generator and one {\it unimodular}
relator, i.e., a relator in which the exponent sum of the new
generator is one:
$$
\^G=\gp{G,t\ |\ w=1}\:=(G*\gp t_\infty)/\!\nc w,
\hbox{ where $w\equiv g_1t^{\epsilon_1}\dots g_nt^{\epsilon_n}$,
\ \
$g_i\in G$,\ \ $\epsilon_i\in\Z$,\ \ and\ \ $\sum\epsilon_i=1$}.
$$

It is known that
a significant part of one-relator group theory extends to
such \emph{unimodular one-relator relative presentations}.
In particular:
\- $G$ embeds (naturally) into $\^G$
[Kl93] (see also [FeR96]);%
\fn{%
However, the natural mapping $G\to\^G$ is never
surjective, except in the case when $w\equiv gt$ [CR01].
}
\- $\^G$ is torsion-free [FoR05];
\- $\^G$ is not simple if it does not coincide with $G$ [Kl05];
\- $\^G$ almost always
(with some known exceptions)
contains
a non-abelian free subgroup
[Kl07];
\- $\^G$ is SQ-universal if $G$
decomposes nontrivially into a free product [Kl06b];
\- the centre of $\^G$ is almost always (with some known exceptions)
trivial
[Kl09].

\enditem
Some generalisations of these results to relative presentations
with several additional generators can be found in
[Kl09],
[Kl07],
[Kl06a],
and [Kl06b].

It is well known that one-relator groups with powered relator are more
similar to free groups than arbitrary one-relator groups. In particular,
Newman's theorem [New68] (see also [LS80]) (reformulated in the modern
language) says that one-relator groups are hyperbolic if the relator is a
proper power. The following recent result is a partial generalisation of
Newman's theorem.

\proclaim{Le Thi Giang's theorem {\rm [Le09]}}.
If a group $G$ is torsion-free,
a word $w\in G*\gp t_\infty$ is unimodular, and
$k\ge2$, then
the group
$$
\~G=\pres<G,t|w^k =1>\:= G*\gp t_\infty/\nc{w^k}
\eqno{(*)}
$$
is
\emph{relatively hyperbolic
{\rm(in the sense of Osin)}
with respect to $G$},
i.e. presentation $(*)$
satisfies a \emph{linear isoperimetric inequality}:
there exists a constant
$C>0$
such that any word $u$ in the alphabet
$G\cup \{t^{\pm1}\}$ representing the identity element of $\~G$
decomposes \(in
$G*\gp t_\infty$\)
into a product of at most $C|u|$
conjugates of $w^{\pm k}$.

\noindent
Henceforth, the symbol
$|u|$ denotes the number of letters $t^{\pm1}$ in the word $u$.

Relatively hyperbolic groups have many good properties.
For example,
they are SQ-universal
(apart from some obvious exceptions) [AMO07],
the word [Far98] and conjugacy [Bum04] problems are solvable
in such groups
(under some natural restrictions).
The same is true for many other algorithmic problems.
More details about relatively hyperbolic groups can be found in
book [Os06].

\smallskip

It turns out that the torsion-freeness condition in
Le Thi Giang's theorem can be replaced by the absence of only order-two
elements. Presently, the following theorem is the unique result about
unimodular relative presentations in which torsion-freeness condition
is weakened to the absence of small-order elements.

\Th.
If
a word $w\in G*\gp t_\infty$ is unimodular and
$k\ge2$, then
the group
$\~G$
defined by relative presentation
$(*)$ contains $G$ as a
(naturally embedded)
subgroup,%
\footnote{**$^)$}{%
In other words, \sl an equation of the form $(w(t))^k=1$,
where $k\ge2$ and the word $w(t)\in G*\gp t_\infty$
is unimodular \emph{is solvable over} any group $G$,
i.e., there exists a group $H$ containing $G$ as a subgroup
and an element $h\in H$ such that $w(h)=1$ in $H$.
}
and
$\gp{G,G^t}=G*G^t$
in~$\~G$.
\newline
If the group $G$ is involution-free or $k\ge3$, then $\~G$ is relatively
hyperbolic with respect to $G$.

\Example 1 {\rm [Le09]}.
The group
$\~G=\pres<g,t|g^3=1, [g,t]^3=1>$ is not hyperbolic
(in particular, it is not relatively hyperbolic with respect to its
finite subgroup $G=\gp g_3$), because the subgroup~$\gp{a^ta,aa^t}$
is a free abelian group of rank two.
This example shows that the unimodularity condition cannot
be omitted from the theorem.

\Example 2.
The Baumslag--Solitar group
$\~G=\pres<g,t|t^g=t^2>$ is not hyperbolic
(in particular it is not relatively hyperbolic with respect to its
cyclic subgroup $G=\gp g$), because the centraliser of the element $t$
is a noncyclic locally cyclic group
$\gp{t^{g^{-1}},t^{g^{-2}},\dots}$.
This example shows that the condition
$k\ge2$ cannot be omitted from the theorem.

\Question.
Can the involution-freeness condition be
omitted from the theorem for
$k=2$?
\newline
\rm
We conjecture that the answer is \emph{no}.

Applying the known facts mentioned above about relatively hyperbolic
groups, we obtain, e.g., the following corollary.

\Corollary 1.
Suppose that a word $w$ is unimodular and either
$k\ge3$ or $k\ge2$ and $G$ is involution-free.  Then
\item{\rm1)}
if $G$ is nontrivial, then $\~G$ is SQ-universal, i.e. any countable group
embeds into a quotient of $\~G$;
\item{\rm2)}
the word and conjugacy problems are solvable in $\~G$ if the corresponding
problems are solvable in $G$ and it is finitely generated.

\Proof
The second assertion follows immediately from the theorem and the results
of Farb [Far98] and Bumagina [Bum04] mentioned above.

To prove the first assertion, it suffices to apply the
Arzhantseva--Minasyan--Osin theorem [AMO07] mentioned above, which says
that a group relatively hyperbolic with respect to its proper subgroup is
either SQ-universal or virtually cyclic.

The group $\~G$ is relatively hyperbolic with respect to $G$ by the
theorem. The subgroup $G\subseteq\~G$ is proper, because
$\~G/\!\nc G=\pres<t|t^k=1>$
is the cyclic group of order $k\ge2$. Finally, $\~G$ is
not virtually cyclic, since according to the theorem $\~G$ contains the
free square of $G$ and it is well-known that the free square of a group of
order larger than two (in particular, any nontrivial group without
involutions) is not virtually cyclic. The remaining case $G\iso\Z_2$
and $k\ge3$ is covered by the Baumslag--Morgan--Shalen theorem [BMS87]
implying that, in this case, $\~G$ contains a non-abelian free subgroup
and, therefore, is not virtually cyclic.

\Corollary 2.
If a word $w$ is unimodular, $G\ne\1$, and $k\ge2$, then $\~G$
contains a non-abelian free subgroup, except in the case where
$G$ consists of two elements, $k=2$, and
$w$ is conjugate in $G*\gp t_\infty$ to a word of the form $gt$, where
$g\in G$ (in this case, $\~G$ is infinite dihedral).

\Proof
According to the theorem, $\~G$ contains the free square of $G$,
which contains a non-abelian free group, except in the case where
$G\iso\Z_2$. In this exceptional case, if $k\ge3$, then the presence
of a non-abelian free subgroup follows from Corollary 1.
If $k=2$, then
the generalized triangle group $\~G=\pres<g,t|g^2=w^2=1>$
satisfies the conditions of a theorem of Howie [How98], which (in
particular) describes generalised triangle groups of such form without
free subgroups.

\Remark.
Our proof shows also that relative presentation $(*)$ is aspherical (if
$w$ is unimodular and $k\ge2$). In particular, this means (see [FoR05])
that each finite subgroup of $\~G$ is conjugate to either a subgroup of
$G$ or a subgroup of the cyclic group $\gp w$.

If we do not assume unimodularity condition in presentation $(*)$
and suppose only that $w$ is not conjugate to elements of $G$ in
$G*\gp t_\infty$, then, as is known, we have, e.g., the following:

\-
the group $G$ embeds naturally
into $\~G$ if either
$G$ is locally indicable [B84],
or $G$ is cyclic and $k\ge2$ [BMS87], [Boy88],
or $k\ge4$ [How90],
or $k\ge3$ and $G$ is involution-free [DuH92];
\-
$\~G$ is relatively hyperbolic with respect to $G$ if
either
$G$ is locally indicable and $k\ge2$ [DuH91]
or $G$ is involution-free and $k\ge4$ [DuH93].

\enditem
A survey of results on one-relator relative presentations with a powered
relator can be found in [DuH93] and [FiR99].

Our approach to the proof of the theorem, as well as
Le Thi Giang's approach, is based on the use of a
standard algebraic trick (Section~2)
and geometric technique:
Howie's diagrams (Section~4)
and car crashes (Sections~6 and~7).
The difference is that we use the crashes
in combination with the weight test, i.e.
the combinatorial Gauss--Bonnet formula (Section~3).
Actually, the major part of the theorem is proven
(in Section 5)
without any ``automobile technique".
The cars are needed only to prove the relative
hyperbolicity when $k=2$ and $G$ is involution-free (Section~8).

\smallskip

\proclaim{Notation}\rm{
which we use is mainly standard.
Note only that if $k\in\Z$, $x$ and $y$ are elements of a
group, and $\phi$ is a homomorphism from this group into another,
then $x^y$, $x^{ky}$, $x^{-y}$, $x^\phi$, $x^{k\phi}$, and
$x^{-\phi}$ denote $y^{-1}xy$, $y^{-1}x^ky$, $y^{-1}x^{-1}y$, $\phi(x)$,
$\phi(x^k)$, and $\phi(x^{-1})$, respectively.
If $X$ is a subset of a group, then
$\gp{X}$ and $\nc{X}$ are the subgroup
generated by $X$ and the
normal subgroup generated by $X$, respectively.
The letters $\Z$, $\N$, and $\R$ denote the set of integers,
positive integers, and real numbers, respectively.
The symbol $\~G$ always denotes the group defined by
presentation~$(*)$}.

The authors thank Le Thi Giang
and an anonymous referee
for useful remarks.

\s2.
An algebraic lemma

The following lemma is an easy generalisation
of Lemma 2.1 from [Le09];
a similar trick with the change of presentation was
used in [Kl93] and later
in many other works
(see, e.g.,
[KP95],
[CG95],
[CG00],
[CR01],
[FeR96],
[FeR98],
[FoR05],
[Kl05],
[Kl06b],
[Kl07],
and
[Kl09]).
A geometric
interpretation of this trick can be found in [FoR05].

\Lemma 1.
If a word $w=g_1t^{\epsilon_1}\dots g_nt^{\epsilon_n}$
is unimodular and cyclically reduced
and $n>1$,
then the group
$\~G$ has a relative presentation of the form
$$
\~G=\pres<H,t|\{p^t=p^\varphi,p\in P\setminus\1\},
\(ct\prod\limits_{i=0}^m(b_i a_i^t)\)^k=1>,
\eqno{(1)}
$$
where $a_i,b_i,c\in H$, $P$ and $P^\phi$ are isomorphic subgroups
of the group $H$, and
$\phi\:P\to P^\phi$ is an isomorphism between them.
In addition,
\item{\rm1)} $m\ge 0$ \(i.e. the product in formula $(1)$ is nonempty\);
\item{\rm2)} $a_i\notin P$ and
             $b_i\notin P^\phi$;
\item{\rm3)} $\gp{P,a_i}=P*\gp{a_i'}$ and
             $\gp{P^\phi,b_i}=P^\phi*\gp{b_i'}$ in $H$,
             where $a_i'\in Pa_i$, $b_i'\in P^\phi b_i$;
\item{\rm4)} the groups $H$, $P$, and $P^\phi$ are free
                products of finitely many isomorphic copies of $G$:
                $H=G^{(0)}*\dots*G^{(s)}$, $P=G^{(0)}*\dots*G^{(s-1)}$,
                and $P^\phi=G^{(1)}*\dots*G^{(s)}$, where $s\ge0$
                                \(if $s=0$, the groups $P$ and $P^\phi$
                are trivial\) and the isomorphism $\phi$ is the shift:
                $\left(G^{(i)}\right)^\phi=G^{(i+1)}$.

\Proof
First, we show that $\~G$ has at least one
presentation of the form~(1) satisfying condition 4).
Since
$\sum\epsilon_i=1$, the word $w$ can be written in the form
$$
w=\left(\prod g_i^{t^{k_i}}\right)t.
$$
Conjugating, if necessary, $w$ by $t$, we can assume that $k_i\ge0$.
Setting $g^{(i)}=g^{t^i}$ for $g\in G$, $G^{(i)}=G^{t^i}$, $s=\max k_i$,
and
$c=\prod g_i^{(k_i)}$, we see that $\~G$ has presentation
$$
\~G\iso\pres<G^{(0)}*\dots*G^{(s)}, t|
\left\{\(g^{(i)}\)^t=g^{(i+1)},\ i=0,\dots,s-1,\ g\in G\right\},\
(ct)^k=1>,
$$
i.e., a presentation of the form (1) (with $m=-1$) satisfying condition 4).

Now, from all presentations of the form (1) satisfying condition 4)
we choose presentations with minimal $s$, and from all these presentations
with minimal $s$ we choose one with minimal $m$. The obtained
presentation (1) is as required.

Indeed, if $m<0$ (i.e., $w=ct$, where $c\in H$), then $s=0$,
because otherwise we might decrease $s$ replacing
all fragments $g^{(s)}$ in the word $c$ by $(g^{(s-1)})^t$. But the
conditions $m<0$ and $s=0$ mean that the initial word $w$ has the form
$w=ct$, where $c\in G$, which contradicts the assumption $n>1$. Thus,
condition 1) holds.

Condition 2) holds because otherwise in presentation (1)
we might replace a fragment $t^{-1}a_it$ with $a_i\in P$
(or a fragment $tb_it^{-1}$ with $b_i\in P^\phi$) by $a_i^\phi$
(or by $b_i^{\phi^{-1}}$, respectively), thereby decreasing $m$ (and not
increasing $s$).

Condition 3) follows from conditions 2) and 4) by virtue of the
following simple fact, whose
proof we leave to the reader as an exercise.

\centerline{\sl
If $u\in A*B$, then $\gp{A,u}=A*\gp{u'}$ for some $u'\in Au$.
}
Lemma 1 is proven.

\Corollary.
If for some $i$
an
equality of the form
$a_i^{n_1}p_1\dots a_i^{n_s}p_s=1$
or
$b_i^{n_1}p_1^\phi\dots b_i^{n_s}p_s^\phi=1$,
where
$s\ge1$,
$n_j\in\Z\setminus\0$,
$p_j\in P$, and
$p_j\ne1$ for $j\ne s$,
holds in $H$,
then
the minimal order of a nonidentity element of $G$
is at most
$\max\limits_{k<l}\left|\sum\limits_{j=k}^l n_j\right|$.

\Proof
This follows immediately from assertions 4), 3), and 2) of Lemma 1.

\s 3.
Maps and weight test

Throughout this paper, the term ``surface" means a closed two-dimensional
oriented surface.

A {\it map} $\Mu$ on a surface $S$ is a finite set of continuous mappings
$\{\mu_i\:D_i\to S\}$, where $D_i$ is a compact oriented two-dimensional
disk, called the $i$th {\it face} or {\it cell} of the map; the boundary
of each face $D_i$ is partitioned into finitely many intervals
$e_{ij}\subset \d D_i$, called the {\it pre-edges} of the map, by a nonempty
set of points $c_{ij}\in \d D_i$, called the {\it corners} of the map.
The images of the corners $\mu_i(c_{ij})$ and the pre-edges $\mu_i(e_{ij})$
are called the {\it vertices} and the {\it edges} of the map, respectively.
It is assumed that
\item{1)}
 the restriction of $\mu_i$ to the interior of each
 face $D_i$ is a homeomorphic embedding preserving orientation; the
 restriction of $\mu_i$ to each pre-edge is a homeomorphic embedding;
\item{2)}
 different edges do not intersect;
\item{3)}
  the images of the interiors of different faces do not intersect;
\item{4)}
 $\bigcup\mu_i(D_i)=S$.

\noindent
Sometimes, we interpret a map $\Mu$ as a continuous mapping
$\Mu\:\coprod D_i\to S$ from a discrete union of disks onto the surface.

The union of all vertices and edges of a map is a graph on the surface,
called the {\it $1$-skeleton}.

We say that a corner $c$ is a corner at a vertex $v$ if $\Mu(c)=v$.  There
is a natural cyclic order on the set of all corners at a vertex $v$; we
call two corners at $v$ {\it adjacent} if they are neighboring with
respect to this order.

By abuse of language, we say that a point or a subset of the surface is
contained in a face $D_i$ if it lies in the image of $\mu_i$. Similarly,
we say that a face $D_i$ is contained in some subset $X\subseteq S$ of
the surface $S$ if $\Mu(D_i)\subseteq X$.

Figure 1 presents a map on the sphere with 10 faces
($A$, $B$, $C$, $D$, $E$, $F$, $G$, $H$, $I$, and $K$),
32 corners,
8 vertices,
16 edges,
and 32 pre-edges.
Note that the number of corners
always equals to the number of pre-edges and is twice the number of edges,
and the value
$$
\chi(S)\:=(the\ number\ of\ vertices)-(the\ number\ of\ edges)+
(the\ number\ of\ faces)
$$
does not depend on the choice of a map on the surface $S$ and is called
{\it the Euler characteristic} of this surface. The Euler characteristic
of the sphere (the only surface of our real interest in this paper)
is two.

\vbox{\vskip-2.5cm
\nobreak

\goodbreak
\bigskip
\centerline{\input 1.PIC}
\nobreak%
\centerline{Fig. \lowercase{1}}%
\goodbreak
\bigskip

}

We need also the following simple but useful fact,
sometimes called the combinatorial Gauss--Bonnet formula.

\proclaim{Weight test \rm[Ger87], [Pri88],
see also [MCW02]}.
If each corner $c$ of
a map on a surface $S$
is assigned a number
$\nu(c)$ {\rm(called the \emph{weight
{\rm or the} value of the
corner
$c$})}, then
$$
\sum_v K(v)+\sum_D K(D)+\sum_e K(e)=2\chi(S).
$$
Here the summations are over all vertices $v$
and all cells $D$ of the map
and the values
$K(v)$, $K(D)$, and $K(e)$, called the \emph{curvatures}
of the corresponding vertex, cell, and edge,
are defined by the formulae
$$
K(v)\:=2-\sum_c \nu(c),
\qquad
K(D)\:=2-\sum_c (1-\nu(c)),
\qquad
K(e)\:=0,
$$
where the first sum is over all corners at the vertex $v$, and
the second sum is over all corners of the cell $D$.

\s 4.
Howie diagrams

Suppose that we have a map $\Mu$ on a surface $S$, the corners of the
map are labeled by elements of a group $H$, and the edges are oriented (in
the figures, we draw arrows on the edges) and labelled by elements of a
set $\{t_1,t_2,\dots\}$ disjoint from the group $H$. The label of a corner
or an edge $x$ is denoted by $\lambda(x)$.

The {\it label of a vertex} $v$ of such a map is defined by the formula
$$
\lambda(v)=\prod_{i=1}^k \lambda(c_i),
$$
where $c_1,\dots,c_k$ are all corners at $v$ listed clockwise.
The label of a vertex is an element of the group $H$ determined up to
conjugacy.
For instance, the label of a vertex in Fig. 1
is
$\lambda(b_2)\lambda(e_1)\lambda(d_1)$.

The {\it label of a face} $D$ is defined by the formula
$$
\lambda(D)=\prod_{i=1}^k
\bigl(\lambda(\Mu(e_i))\bigr)^{\epsilon_i}\lambda(c_i),
$$
where $e_1,\dots,e_k$ and $c_1,\dots,c_k$ are all pre-edges and all
corners of $D$ listed anticlockwise, the endpoints of $e_i$ are
$c_{i-1}$ and $c_i$ (subscripts are
modulo $k$), and $\epsilon_i=\pm1$ depending on whether the homeomorphism
$e_i\mathop\to\limits^\Mu\Mu(e_i)$ preserves or reverses orientation.
Simply speaking, to obtain the label of a face, we should go around its
boundary anticlockwise, writing out the labels of all corners and edges we
meet; the label of an edge traversed against the arrow should be raised to
the power $-1$.

The label of a face is an element of the group $H*F(t_1,t_2,\dots)$ (the
free product of $H$ and the free group with basis $\{t_1,t_2,\dots\}$)
determined up to a cyclic permutation. More precisely, the right-hand side
of our formula for $\lambda(D)$ is called the {\it label of the face $D$
written starting with the pre-edge $e_1$}.

For instance, if the label of each edge in Fig. 1 is $t$, then the label
of the face $B$ written starting with the pre-edge~$\alpha$ is
$
t\lambda(b_1)t\lambda(b_2)t^{-1}\lambda(b_3).
$

Such a labelled map is called a {\it Howie diagram} (or
simply {\it diagram}) over a relative presentation
$$
K=\gp{H,t_1,t_2,\dots\ |\ w_1=1,w_2=1,\dots}
\eqno{(**)}
$$
if
\item{1)}
  some vertices and faces are distinguished and called {\it
  exterior}; the remaining vertices and faces are called {\it interior};
\item{2)}
    the label of each interior face is a cyclic permutation of one of
  the words $w_i^{\pm1}$;
\item{3)}
  the label of each interior vertex is the identity element of $H$.

\enditem
Figure 4 presents all possible interior faces of Howie diagrams
over presentation (1).

A diagram is said to be {\it reduced} if it contains no such
edge $e$ that both faces containing $e$ are interior, these faces
are different and the label
of one of these face written
 starting with the label of
$e$ is inverse
to the label of the other face written ending with
the label of $e$;
such a pair of faces with a common edge is
called a {\it reducible pair}. For example, the faces $D$ and $E$ in Fig.1
form a reducible pair if
$\lambda(d_i)=(\lambda(e_i))^{-1}$
and the labels of all
edges are equal.

The following lemma is an analogue of the van Kampen lemma for relative
presentations.

\Lemma 2
{\rm[How83]}.
The natural mapping from a group $H$ to the group with relative
presentation $(**)$ is noninjective if and only if there exists a
spherical diagram over this presentation with no exterior faces and a
single exterior vertex whose label is not 1 in $G$. A
minimal \(with respect to the number of faces\) such diagram is
reduced. \hfil\break
If this natural mapping is injective, then we have
the equivalence: the image of an element $u\in
H*F(t_1,t_2,\dots)\setminus \1$ is 1 in the group $(**)$ if and only
if there exists a spherical diagram over this presentation without
exterior vertices and with a single exterior face with label $u$. A
minimal \(with respect to the number of faces\) such diagram is
also reduced.

Diagrams on the sphere with a single exterior face and no exterior
vertices are also called {\it disk diagrams}, the boundary of the exterior
face of such a diagram is called the {\it contour} of the diagram.

Let $\phi\:P\to P^\phi$ be an isomorphism between two subgroups of a group
$H$.  A relative presentation of the form
$$
\gp{H,t\ |\ \{p^t=p^\phi;\
p\in P\setminus\1\}, w_1=1,\ w_2=1,\ \dots}
\eqno{({**}*)}
$$
is called a {\it
$\phi$-presentation}. A diagram over a $\phi$-presentation $({**}*)$ is
called {\it $\phi$-reduced} if it is reduced and different interior cells
with labels of the form $p^tp^{-\phi}$, where $p\in P$, have no common
edges.

\Lemma 3 {\rm[Kl05]}. A minimal \(with respect to the number of faces\)
diagram among all spherical diagrams over a given $\phi$-presentation
without exterior faces and with a single exterior vertex with nontrivial
label is $\phi$-reduced. If no such diagrams exists, then a minimal
diagram among all disk diagrams with a given label of contour is
$\phi$-reduced.  \rm In other words, the complete $\phi$-analogue of Lemma
2 is valid.

The idea of the proof is shown in Fig.2.

\vfil\break

\

\vskip-10mm

\goodbreak
\bigskip
\centerline{\input 2.PIC}
\nobreak%
\centerline{Fig. \lowercase{2}}%
\goodbreak
\bigskip

A relative presentation ($\phi$-presentation) over which there exists
no reduced ($\phi$-reduced) spherical diagrams with
no exterior faces and
a single
exterior vertex are called {\it aspherical} (respectively, {\it
$\phi$-aspherical}).

Suppose that we have a map on a surface all whose edges are oriented
(e.g., a Howie diagram). Such a map has 4 kinds of corners:
$(++)$, $(--)$, $(+-)$, and $(-+)$ (Fig. 3).

\goodbreak
\bigskip
\centerline{\input 3.PIC}
\nobreak%
\centerline{Fig. \lowercase{3}}%
\goodbreak
\bigskip

The following lemma is obvious.

\Lemma 4. In the anticlockwise listing of the corners at a vertex
$v$, the corners of type $(++)$ alternate with corners of type $(--)$. If
at a vertex $v$ there are no corners of type $(++)$, or, equivalently,
there are no corners of type $(--)$, then either all corners at $v$ are of
type $(+-)$ \(in this case, $v$ is called a {\it sink}\), or all corners at
$v$ are of type $(-+)$ \(in this case, $v$ is called a {\it source}\).

\vskip-5mm

\goodbreak
\bigskip
\centerline{\input 4A.PIC}
\nobreak%
\centerline{Fig. \lowercase{4A}}%
\goodbreak
\bigskip

\vskip-1mm

\goodbreak
\bigskip
\centerline{\input 4B.PIC}
\nobreak%
\centerline{Fig. \lowercase{4B}}%
\goodbreak
\bigskip

\vskip-1mm

\goodbreak
\bigskip
\centerline{\input 4C.PIC}
\nobreak%
\centerline{Fig. \lowercase{4C}}%
\goodbreak
\bigskip

\s 5.
The proof of a major part of the theorem

In this section, we prove all assertions of the theorem except the
relative hyperbolicity for $k=2$.

If the word $w$ is conjugate to a word $gt$, then the group
$\~G$ is the free product of the group
$G$ and a cyclic group of order $k$, and
all assertions of the theorem are obvious.
If the letters $t^{\pm1}$ occur more than once in the word $w$,
then, by Lemma~1, the group $\~G$ has presentation (1).

Consider a
$\phi$-reduced
spherical Howie diagram over presentation (1) that has either no exterior
faces and one exterior vertex or
no exterior vertices and one exterior face.
Faces with label of the form $p^{-\varphi}p^t$ are called \emph{digons},
the other interior faces are called \emph{large faces}.

Vertices and edges belonging to the boundary of the exterior face are
called \emph{boundary}. The exterior vertex (if it exists) is also
considered as a boundary vertex.

A digon is called \emph{special} if
its both neighboring faces are interior and
one of its corners (called \emph{positive})
is adjacent with corners of types $(++)$ and $(--)$
(Fig.5).
Note that the other corner of a special digon
(called \emph{negative}) is automatically
non-adjacent with corners of type $(++)$ and $(--)$.

\goodbreak
\bigskip
\centerline{\input 5.PIC}
\nobreak%
\centerline{Fig. \lowercase{5}}%
\goodbreak
\bigskip

Let us assign a value (weight)
$\nu(\gamma)$ to each corner $\gamma$ of the diagram by the following
rule:
$$
\nu(\gamma)=\cases{
0 &if $\gamma$ is a corner of a nonspecial digon
\cr
&or a corner of type $(++)$ or $(--)$ of an
interior face (the label of such a corner is $c^{\pm1}$);
\cr
-1 &if $\gamma$ is a negative corner of a special digon;
\cr
1, &otherwise.
\cr
}
$$

Let us calculate the curvatures of vertices and faces according to the
weight test (see Section 3). For faces, we have
$$
K(\hbox{digon})=0,
\quad
K(\hbox{large face})=2-k,
\quad
K(\hbox{exterior face})=2.
$$
For a vertex $v$, the curvature is
$$
K(v)=2
+n-l-
{\bf p}
-x,
\eqno{(2)}
$$
where $l$ is the number of corners of types $(+-)$ and $(-+)$ of
large faces,
${\bf p}$
is the number of positive corners of special digons,
$n$ is the number of negative corners of special digons,
and $x$ is the number of corners of the exterior face
(all corners are at the vertex $v$).

Each negative corner of a special digon is adjacent to
two corners of type $(+-)$ or $(-+)$ of large faces
(by the definition of special digons),
and no corner of type $(+-)$ or $(-+)$ can be adjacent to two
negative corners
(since otherwise, the corresponding large face would have both
a corner of type $(++)$ and a corner of type $(--)$).
Therefore,
$l\ge 2n$.

Note also that corners of types $(++)$ and $(--)$ at
a non-boundary
vertex
alternate (Lemma 4) and
cannot be adjacent
(since the diagram is reduced):
between two such corners there must be a corner of weight 1
(either a corner of type $(+-)$ or $(-+)$ of a large face or a
positive corner of a special digon).
Taking into account the preceding remark about negative corners,
we conclude that the sum of weights of corners lying between corners of
type $(++)$ and $(--)$ (if we list them clockwise around the vertex $v$)
is at least
one (Fig.~6, left). Therefore,
a non-boundary
vertex with positive
curvature must be either a source or a sink and, for such vertex,
${\bf p}=0$,
and
either $n=1$ and $l=2$ or $n=0$ and $l=1$ or $n=0$ and $l=0$
($n<2$, since otherwise, formula (2) and the inequality
$l\ge 2n$ mentioned above would give a nonpositive curvature).
See Fig.~6, the boldface digits denote the values of corners.

\goodbreak
\bigskip
\centerline{\input 6.PIC}
\nobreak%
\centerline{Fig. \lowercase{6}}%
\goodbreak
\bigskip

The first case
($n=1$ and $l=2$)
for
a non-boundary
vertex is impossible, because
the label of such a vertex, i.e., the product of labels of
corners, is
$a_m^{-1}p_1a_mp_2$ (if the vertex is a source) or
$b_0^{-1}p_1^\phi b_0p_2^\phi$ (if the vertex is a sink), where $p_1$ and
$p_2$ lie in $P$ and are not 1 (since the diagram is reduced)
and, therefore, the label of the vertex is not 1 by Corollary
of Lemma 1; thus this vertex can not be interior. The second and third
cases
($n=0$ and $l\in\{0,1\}$)
for
a non-boundary
vertex are impossible by nearly the same reason:
they would imply an equality of the form $a_i^{\pm1}p_1=1$,
$b_i^{\pm1}p_1^\phi=1$,
$p_2=1$, or $p_2^\phi=1$, where $p_1\in P\ni p_2\ne1$.

Thus, the curvature of any
non-boundary
vertex $v$ is
nonpositive. The curvatures of interior faces are
also nonpositive (for $k\ge2$), the curvature of a boundary vertex is at
most two (this follows from formula (2) and the inequality $l\ge 2n$),
while the
total curvature must be four according to the weight test.

This means that, first, there exist no diagrams
without exterior faces and with single exterior vertex,
i.e. the
natural mapping $H\to\~G$
(and, hence, the natural mapping $G\to\~G$) is
injective by Lemma~2; and secondly, if there is one exterior face and no
exterior vertices and $k\ge3$, then the
number of interior large faces
is bounded by a linear function of the perimeter of the exterior face:
$$
2\cdot(\hbox{the perimeter of the exterior face})-
(k-2)\cdot(\hbox{the number of large interior faces})+2\ge4.
$$

It is easy to see that
such an isoperimetric inequality for presentation~(1) implies
the usual linear
isoperimetric inequality for presentation $(*)$ (see [Le09]), i.e., the
relative hyperbolicity of $\~G$ for $k\ge3$.
For the sake of completeness, we prove this fact here.
\Proposition 1.
Suppose that some word $u\in G*\gp t_\infty$
represents the identity element of the group $\~G$, i.e. $u$
can be represented as a product of the form
$$
u=v_1\dots v_p w_1\dots w_s,
$$
where each $v_i$ is conjugate to
a word of the form $p^{-t}p^\phi$ ($p\in P$)
in the group
$H*\gp t_\infty$ and
each $w_i$ is conjugate to the word
$\(ct\prod\limits_{i=0}^m(b_i a_i^t)\)^{\pm k}$
in
$H*\gp t_\infty$
(in the notation of
Lemma 1, where $G$ is the same as $G^{(0)}$).
Then $u$ can be represented as a product of $s$
words conjugate
to $w^{\pm k}$
in $G*\gp t_\infty$.
\newline
\rm
Informally, any isoperimetric inequality for presentation
(1) counting only long relators
(only large faces) implies the same isoperimetric inequality
for presentation $(*)$.

\Proof
In the group $\pres<H,t|\{p^t=p^\phi\;;\; p\in P\}>$
(isomorphic to $G*\gp t_\infty$), the words
$v_i$ represent the identity element ant the words $w_i$ are conjugate to
$w^{\pm k}$ (because $ct\prod\limits_{i=0}^m(b_i a_i^t)$ is equal to
a cyclic shift of $w$ by the construction). This implies the assertion of
Proposition 1.

\medskip

Resuming the proof of the theorem, let us show
that $\gp{G,G^t}=G*G^t$ in the group $\~G$.
If $H\ne G$, i.e., if $P\ne\1$, i.e., if $s>0$ in Lemma 1, then we have
nothing to prove, because it is already proven that the natural mapping
$H=G*G^t*\dots\to\~G$ is injective.

It remains to consider the
case $H=G$ (i.e., $P=\1$).
Suppose that $u\in G*G^t$ is a reduced nonempty word
representing the identity element of $\~G$. By Lemma~2, $u$ is the
label of the exterior face of some $\phi$-reduced spherical diagram over
presentation $(*)$ (which coincides with presentation (1) in the case
under consideration) without exterior vertices and with a unique exterior
face. Since digons are absent and the exterior face has no corners of
types $(++)$ and $(--)$, the curvature of each boundary vertex is
nonpositive. The sum of curvatures of all faces and vertices must be four,
but the unique positive term in this sum is two (the curvature of
the exterior face). This contradiction with the weight test completes
the proof of the theorem, except the assertion about relative
hyperbolicity for $k=2$.

\Remark.
This argument proves also the $\phi$-asphericity
of presentation (1) (for $k\ge2$), which implies (see [FoR05])
the
asphericity of
presentation~$(*)$.

In the remaining part of this paper, we prove relative
hyperbolicity for $k=2$.

\s 6.
Motions

All definitions and facts of this section
are taken
from paper [Kl05].

Consider a map $\Mu$ on a closed oriented surface $S$.
Some corners
of this map are distinguished and called
\emph{stop corners}.

A \emph{car} moving around a face $D$ of this map is a
continuous locally nondecreasing
\fn{%
We call a continuous mapping $\alpha\:X\to Y$ from an oriented
circle $X$ to an oriented circle $Y$ (locally) {\it nondecreasing} if
the preimage of any interval $U\subset Y$ is a union of intervals such
that the restriction of $\alpha$ to each of these intervals is
a nondecreasing function (in the usual sense, as a function from one
oriented interval to another).
}
mapping from an oriented circle
$R$ (\emph{the circle of time})
to the boundary $\d D$ of the face $D$ such that the
preimage of each point, except possibly stop corners,
is discrete.

Simply speaking, each car moves without U-turns
and infinite decelerations and accelerations along the
boundary of its face anticlockwise, possibly stopping for a finite time at
some corners.
And this motion is periodic.

We say that a car $\alpha_i$ {\it is at} a corner $c\in\d D_i$ at a moment
of time $t\in R$ if $\alpha_i(t)=c$; we also say that a car
$\alpha_i$ {\it is at} a point $p\in S$ at a moment $t\in R$ if
$\mu_i(\alpha_i(t))=p$. If the number of cars being at a moment $t\in R$ at
a point~$p$ of the 1-skeleton of $S$ equals the multiplicity of this
point (in other words, $\bigcup\alpha_i(t)\supseteq\Mu^{-1}(p)$), then we
say that at the point~$p$ at the moment $t$ a {\it complete collision}
occurs; the point $p$ is called a {\it point of complete collision}.
Points of complete collision lying on edges are called simply {\it points
of collision}.

\emph{A multiple motion of period $T$ with separated stops}
on a map $\Mu$
is a
set of cars
$\alpha_{D,j}\:R\to\d D$, where $j=1,\dots,d_D$,
such that
\item{1)}
$d_D\ge 1$ (i.e. each face is moved around by at least one
car);
\item{2)} at each vertex $v$ at which there are stop corners,
          the stops are separated in the following sense: let
          $c_1,\dots,c_k$ be all stop corners at $v$ enumerated
          anticlockwise; it is required that, for each $i$, at corners
          $c_{i}$ and $c_{i+1}$ (subscripts are modulo~$k$), cars are
          never located simultaneously. (In particular, this implies that
          $k\ge 2$.)
\item{3)}
$\alpha_{D,j}(t+T)=\alpha_{D,j+1}(t)$ for any $t\in R$ and
$j=\{1,...,d_D\}$ (subscripts are modulo ${d_D}$,
and the addition of points of the circle~$R$ is defined naturally:
$R=\R/l\Z$);

\item{4)}
there exists a partition of each circle $\partial D$
into $d_D$ arcs
(with disjoint interiors)
such that during the interval of time $[0,T]$ each
car $\alpha_{D,j}$ moves along the $j$th arc.

\proclaim{Car-crash test \rm [Kl05], [Kl97]}.
For any multiple motion with separated stops
on a map $\Mu$ on a
surface $S$, we have
$$
\sum_v K'(v)+\sum_e K'(e)+\sum_D K'(D)=\chi(S),
$$
where
the sums are over all vertices $v$, edges $e$, and
faces $D$ of the map $\Mu$.
\newline
Here
$K'(D)=1-d_D$,
the value
$K'(e)$ is the number of collision points on an
edge $e$
(not counting the end-points),
and
$K'(v)=1$ if at the vertex $v$ a complete collision occurs; otherwise
$K'(v)$ is an integer nonpositive number
{\rm (whose exact definition
can be found in [Kl05])}.

Throughout this paper,
the surface is always the sphere, its Euler characteristic is 2.

\s 7.
Standard multiple motion

In this section, we define some particular
multiple motion on Howie diagrams
over presentation (1). Our definition almost
literally
repeats a definition from [Le09]. A similar motion
was considered in [Kl05].

The following motion on a Howie diagram over
presentation (1) is called \emph{standard}:

\item{a)}
  the car going around an interior face with label $p^{-\varphi}p^t$ moves
  anticlockwise
  uniformly with unit speed (one edge per a unit time)
  visiting the corner of type $(+-)$ at the even moments of time
  (Fig. 4a);

\item{b)}
  An interior face with label
  $\(ct\prod\limits_{i=0}^m b_ia_i^t\)^k$
  are moved around by $k$ cars;
  for $m>0$, they stay at
  the corners of type $(++)$ during the time intervals
  $[2m+2,4m+1]+(4m+2)\Z$, and moves anticlockwise uniformly with unit speed
  all the remaining time; for $m=0$, each car moves without stops with
  speed 2 when it moves in the direction of an edge,
  and with speed 1 when it moves against the direction of an edge;
  at time zero the car is at
  a corner of type $(+-)$ (Fig. 4b);

\item{c)}
  An interior face with label
  $\(ct\prod\limits_{i=0}^m b_ia_i^t\)^{-k}$
  are moved around by $k$ cars;
  for $m>0$, they stay at
  the corners of type $(--)$ during the time intervals
  $[1,2m]+(4m+2)\Z$, and moves anticlockwise uniformly with unit speed
  all the remaining time; for $m=0$, each car moves without stops with
  speed 2 when it moves against the direction of an edge,
  and with speed 1 when it moves in the direction of an edge;
  at time zero the car is at
  a corner of type $(+-)$ (Fig. 4c);

\item{d)}
  An exterior face is moved around by one car;
  it moves with period $4m+2$; at time zero, it
  is at some vertex; during the interval
  $[0, {1\over4}]$, it (rapidly) moves counterclockwise along the entire
  boundary of the face, except the last edge; and at the remaining time
  it (slowly) goes along this edge.

\enditem
The standard motion is periodic with period $4m+2$
(on faces with label $p^{-\varphi}p^t$ minimal period is two).
Figure~4
shows the detailed schedule of the motion of
cars moving around interior cells
during the
interval $[0,4m+2)$; the framed numbers near edges denote
the speed of
the cars on these edges (the default speed is unit).

\Lemma 5 {\rm(cf. [Le09], [Kl05])}.
Suppose that a Howie diagram over presentation (1)
has at most one
exterior face. Then
the standard motion is
a motion with separated stops. Complete collisions which
occur not on the boundary of the exterior face can occur only at
vertices being sinks or sources and only at integer moments of
time.
On each edge of the boundary of the exterior face there are at most
$k(2m+1)$ points of complete collision.

\Proof
Let us declare all corners of types $(++)$ and $(--)$ to be stop
corners.
The schedule of the standard motion is such that cars are never
located simultaneously at corners of types $(++)$ and $(--)$: the
corners of type $(--)$ are visited only during the first half of
the period, while
the corners of type $(++)$ are visited during the second half of the period.
The car moving around the exterior face is not at corners at all
at such moments.
This and Lemma 4 imply that the standard motion is a motion with
separated stops.
A collision on an edge separating two interior faces
at a moment $t$ means that at this moment
the direction of the motion of one of the cars coincides with the
direction of the edge, while the direction of the motion of the other
colliding car is opposite to the direction of the edge.
But the schedule of the standard motion is such that, at each moment $t$,
either all cars
moving around interior faces and
being on edges move in the direction of the edge (this is
so when the integer part of $t$ is odd), or all cars being on edges move
in the direction opposite to the direction of the edge (this is so when
the integer part of $t$ is even).
Note also that the definition of multiple motion implies that there are no
overtakings.
Therefore, collisions can occur only at
vertices; the separatedness of stops implies that a vertex of complete
collision can not have stop corners and, therefore, is a source or a sink.
The cars visit such vertices only at integer moments of time (even for
sinks and odd for sources).

The car $\beta$ moving around the exterior face can collide with at most
$k$ cars
on each edge $e$.
During the period
$[0;4m+2)$
the car $\beta$
occurs
on each edge
only once,
while each car moving along this edge in the
opposite direction
occurs on $e$ at most $2m+1$ times
(this value is attained on
digons). Therefore, during the period, on each edge of the boundary
of the exterior face at most $k(2m+1)$ collisions occur.
This very rough estimate completes the proof.

\s 8.
Completion of the proof of the theorem

In this section, we complete the proof of the theorem, i.e., we prove
that $\~G$ is relatively hyperbolic with respect to
$G$ if $G$ contains no involutions and $k=2$ (however, the proof below
is suitable for any $k\ge2$).

If the word $w$ is conjugate to $gt$, then
$\~G$ is the free product of
$G$ and the cyclic group of order $k$, and
we have nothing to prove.
If letters $t^{\pm1}$ occur more then twice in $w$, then by Lemma~1
$\~G$ has presentation (1).

Consider a
$\phi$-reduced
spherical Howie diagram over presentation (1)
without exterior vertices and with one exterior face.
As in Section~5, it suffices to show that the diagram satisfies
a linear isoperimetric inequality,
i.e., the number of large interior faces is bounded by a linear
function of the perimeter of the exterior face.

Let us assign a value (weight) to each corner of the diagram
as in Section 5. Recall that, for such weights, the curvatures of
interior vertices are nonpositive.
Moreover, according to formula (2),
the curvature of an interior vertex
can be zero only in the following
cases:
\item{a)} $p>0$ (and, therefore, the vertex is neither a source nor
a sink);
\item{b)} $p=0$, $n=0$, $l=2$;
\item{c)} $p=0$, $n=1$, $l=3$;
\item{d)} $p=0$, $n=2$, $l=4$
(Fig. 7).

\vskip-13mm

\hbox{\qquad\qquad\vbox{
\goodbreak
\bigskip
\centerline{\input 7.PIC}
\nobreak%
\centerline{Fig. \lowercase{7}}%
\goodbreak
\bigskip
}}

Note that in cases a), b), c), d)
a complete collision cannot occur
at the vertex
$v$
under the standard motion (Section~7).
Indeed, by virtue of Lemma~5,
a vertex of complete collision must be either a source or a sink;
therefore, in case a) we have no complete collision.
A complete collision in case b), when the vertex $v$
is, e.g., a source, would imply, according to the schedule of the
motion, that both corners of large faces at this vertex have labels
$a_i^{\pm1}$ with the same subscript~$i$,
and the product of all these labels is 1 in the group $G$; this is
impossible
by virtue of the reducedness of the diagram, the absence of involutions, and
the
corollary of Lemma~1. For the same reason, complete collisions cannot
occur in cases c) and d): in these cases, all corners of large faces must
have labels $a_m^{\pm1}$ if the vertex is a source, or $b_0^{\pm1}$ if
the vertex is a sink.

Note also that, for the standard motion (Section 7),
we have
$$
K'(\hbox{digon})=0,
\
K'(\hbox{large face})=1-k,
\
K'(\hbox{non-boundary edge})=0,
\
K'(\hbox{boundary edge})\le k(2m+1),
$$
where the value $K'$ is defined in Section~4 (car-crash test).
The last two inequalities follow from Lemma~5.

Now, we define the
\emph{combined curvature} of vertices, faces, and edges by the formula
$$
K_\Sigma(\cdot)\:=K(\cdot)+K'(\cdot).
$$
Clearly,
$
K_\Sigma(v)\le0
$
for any interior vertex $v$, because
$K(v)$ is either a negative integer or zero, but in the latter
case, as we have seen, there are no complete collision at the vertex $v$
and, therefore, $K'(v)\le0$.

It remains to note that, for any non-boundary edge $e$ and
any interior large face~$\Gamma$,
$$
K_\Sigma(e)=K'(e)=0
\ \hbox{and}
\quad
K_\Sigma(\Gamma)=K(\Gamma)+K'(\Gamma)=
2-k+(1-k)\le-1
\ \hbox{for $k\ge2$}.
$$
The combined curvature of a boundary edge is bounded by some constant
(depending only on $k$ and $m$) by Lemma~5. The combined curvature
of a
boundary vertex is at most three
(since $K(v)\le2$, as mentioned in Section~5). The combined
curvature of the exterior face is two.
On the other hand, the sum of the combined curvatures of all vertices,
edges, and faces must be
$4+2$, according to the weight test and the car-crash test.

This means that
the number of interior large faces
is bounded by a linear function of the perimeter of the exterior face:
$$
(D+3)\cdot(\hbox{perimeter of the exterior face})-
(\hbox{number of large interior faces})+2\ge 4+2,
$$
where $D=k(2m+1)$ is the constant from Lemma~5 (this is a very
rough estimate). This isoperimetric inequality completes the proof
(by virtue of Proposition 1).

Other applications of the combined test and a description of all possible
tests (in some exact sense) can be found in [Kl97].

\REFERENCES

\[B84]
Brodskii S.D.
Equations over groups and one-relator groups
{// Sib. Mat. Zh.} 1984. {T.25}. no.2. P.84--103.

\[Kl05]
Klyachko Ant. A.
The Kervaire--Laudenbach conjecture and presentations of simple groups
{// Algebra i Logika}. 2005. {T. 44}. {no.4}. P. 399--437.
See also
{arXiv:math.GR/0409146}

\[Kl06a]
Klyachko Ant. A.
How to generalize known results on equations over groups
{// Mat. Zametki}. 2006. {T.79}. no.3. P.409--419.
See also
arXiv:math.GR/0406382.

\[Kl06b]
Klyachko Ant. A.
The SQ-universality of one-relator relative presentations
{// Mat. Sbornik}. 2006. {T.197}. no.10. P.87--108.
See also
arXiv:math.GR/0603468.

\[Kl07]
Klyachko Ant. A.
Free subgroups of one-relator relative presentations
{// Algebra i Logika}. 2007. {V.46}. no.3. P.290--298
See also
arXiv:math.GR/0510582.

\[KP95]
Klyachko Ant. A., Prishchepov M. I.
The descent method for equations over  groups
{// Moscow Univ. Math. Bull.} 1995, {V.50}  P. 56--58.

\[LS77]
Lyndon R.C., Schupp P.E.
{Combinatorial Group Theory},
Springer-Verlag, Berlin/Heidelberg/New~York, 1977.

\[AMO07]
Arzhantseva G., Minasyan A., Osin D.
The SQ-universality and residual properties of relatively hyperbolic
groups
//
Journal of Algebra. 2007. 315:1, 165--177.
See also
arXiv:math.GR/0601590

\[BMS87]
Baumslag G., Morgan J.W., Shalen P.B.
Generalized triangle groups
//
Math. Proc. Camb. Phil. Soc. 1987. 102. 25-31.


\[Boy88]
Boyer S.
On proper powers in free products and Dehn surgery
//
J. Pure Appl. Algebra. 1988. 51:3. 217--229.

\[Bum04]
Bumagina I.
The conjugacy problem for relatively
hyperbolic groups
//
Algebraic \& Geometric Topology.
2004. V.4. P.1013--1040.
See also
arXiv:math/0308171

\[CG95]
Clifford A., Goldstein R.Z.
Tesselations of $S^2$ and equations over torsion-free groups
{// Proc. Edinburgh Math. Soc.} 1995. {V.38}. P.485--493.

\[CG00]
Clifford A., Goldstein R.Z.
Equations with torsion-free coefficients
{// Proc. Edinburgh Math. Soc.} 2000. {V.43}. P.295--307.

\[CR01]
Cohen M. M., Rourke C.
The surjectivity problem for one-generator, one-relator extensions of
torsion-free groups
{// Geometry \& Topology}. 2001. {V.5}. P.127--142.
See also
arXiv:math.GR/0009101

\[DuH91]
Duncan A.J, Howie J.
The genus problem for one-relator products of locally indicable groups//
Mathematische Zeitschrift. 1991. 208:1. 225--237

\[DuH92]
Duncan A.J, Howie J.
Weinbaum's conjecture on unique subwords of nonperiodic words
//
Proc. Amer. Math. Soc. 1992. 115. 947-954.

\[DuH93]
Duncan A.J, Howie J.
One-relator products with high-powered relator,
in: Geometric group theory
(G.A.Niblo, M.A.Roller, eds.),
P.48--74,
Cambridge Univ. Press, Cambridge (1993).

\[Far98]
Farb B.
Relatively hyperbolic groups
// GAFA. 1998. V.8. 810--840.

\[FeR96]
Fenn R., Rourke C.
Klyachko's methods and the solution of equations over torsion-free groups
{// L'Enseignment Math\'ematique.} 1996. {T.42}. P.49--74.

\[FeR98]
Fenn R., Rourke C.
Characterisation of a class of equations with solution over torsion-free
groups,
from {``The Epstein Birthday Schrift"},
{(I. Rivin, C. Rourke and C. Series, editors)},
{Geometry and Topology Monographs.} 1998. {V.1}. P.159-166.

\[FiR99]
Fine B., Rosenberger G.
Algebraic generalizations of discrete groups. A path to combinatorial
group theory through one-relator products.
Monographs and Textbooks in Pure and Applied
Math. 223. Marcel Dekker, Inc., New York, 1999. x+317 pp.

\[FoR05]
Forester M., Rourke C.
Diagrams and the second homotopy group
{// Comm. Anal. Geom.} 2005. {V.13}. P.801-820.
See also
arXiv:math.AT/0306088

\[Ger87]
Gersten S.M.
Reducible diagrams and equations over groups.
In {Essays in group theory}, P.15--73.
Springer, New York-Berlin, 1987.

\[How83]
Howie J.
The solution of length three equations over groups
{// Proc. Edinburgh Math. Soc.} 1983. {V.26}. P.89--96.

\[How90]
Howie J.
The quotient of a free product of groups by a single high-powered relator.
II. Fourth powers
//
Proc. London Math. Soc. 1990. 61. 33--62.

\[How98]
Howie J.
Free subgroups in groups of small deficiency
{// J. Group Theory}. 1998. V.1. no. 1.  P.95--112.

\[Kl93]
Klyachko Ant. A.
A funny property of a sphere and equations over groups
{// Comm. Algebra}. 1993. {V.21}. P.2555--2575.

\[Kl97]
Klyachko Ant. A.
Asphericity tests
{// IJAC}. 1997. {V.7}. P.415--431.

\[Kl09]
Klyachko Ant. A.
The structure of one-relator relative presentations and their centres
//Journal of Group Theory, 2009, 12:6, 923--947.
See also
arXiv:math.GR/0701308

\[Le09]
Le Thi Giang.
The relative hyperbolicity of one-relator relative presentations
//
Journal of Group Theory. 2009. 12:6, 949--959.
See also
arXiv:0807.2487

\[MCW02]
McCammond J.P., Wise D.T.
Fans and ladders in small cancellation theory.
{Proc. London Math. Soc. (3)}. 2002. 84(3):599--644.

\[New68]
Newman B.B.
Some results on one-relator groups
{// Bull. Amer. Math. Soc.} 1968. {V.74}. P.568--571.

\[Os06]
Osin D.V.
Relatively hyperbolic groups:
Intrinsic geometry, algebraic properties, and algorithmic problems.
Memoirs Amer. Math. Soc. 179 (2006), no. 843, vi+100 pp.
See also
arXiv:math/0404040

\[Pri88]
Pride S.J.
Star-complexes, and the dependence problems for hyperbolic complexes.
{Glasgow Math. J.} 1988. 30(2):155--170.

\end